\documentclass{book}
\usepackage{amsmath, amssymb, tikz-cd, hyperref, url}
\usepackage[affil-it]{authblk}
\title{A Very Short Introduction to Topos Theory\\ (adapted from Prof. Pettigrew's notes)}
\author{Eric Schmid}
\affil{New York University}
\date{}

\begin{document}

\maketitle
\frontmatter
\chapter{Preface}
The original notes by Prof. Pettigrew are invaluable, and Prof. Awodey's book is indispensible. My deepest gratitude to Prof. Pettigrew for making his notes publicly available online and massive thanks to Prof. Awodey for his rigorous work in category theory. I am extremely grateful for the academic support from Prof William Chin. 

\tableofcontents
\mainmatter
\chapter{Motivating Category Theory}
\section{The Idea Behind Category Theory}
Category theory encourages a shift from focusing on the intrinsic properties of mathematical objects to emphasizing their roles and interactions. This approach can be summarized by the phrase: "Ask not what a thing is; ask what it does" \cite{Pettigrew}. This shift in perspective helps unify various areas of mathematics by providing a common framework to describe different mathematical constructs through their relationships.

\section{Traditional vs. Categorical Perspectives}
Traditional mathematics often describes objects by their intrinsic properties:
\begin{itemize}
    \item \textbf{Groups:} Defined by a set with an operation following group axioms.
    \item \textbf{Sets:} Collections of distinct elements.
    \item \textbf{Topological Spaces:} Sets endowed with a topology specifying open sets.
\end{itemize}

In contrast, category theory interprets these objects through their morphisms:
\begin{itemize}
    \item \textbf{Groups:} Understood through group homomorphisms.
    \item \textbf{Sets:} Understood through functions.
    \item \textbf{Topological Spaces:} Understood through continuous maps.
\end{itemize}

Category theory provides a more flexible and general approach by focusing on the relationships between objects rather than their internal structure.

\chapter{The Definition of a Category}
A category $\mathcal{C}$ consists of:
\begin{itemize}
    \item A collection of objects $\mathrm{Ob}(\mathcal{C})$.
    \item A collection of morphisms (arrows) $\mathrm{Ar}(\mathcal{C})$.
    \item Two functions: domain and codomain, assigning to each morphism $f: A \rightarrow B$ its source $A$ and target $B$.
    \item A composition function $\circ$ that assigns to each pair of composable morphisms $f: A \rightarrow B$ and $g: B \rightarrow C$ an arrow $g \circ f: A \rightarrow C$.
\end{itemize}

\[
\begin{tikzcd}
A \arrow[r, "f"] & B \arrow[r, "g"] & C
\end{tikzcd}
\]

These components must satisfy the following axioms:
\begin{enumerate}
    \item \textbf{Associativity:} For all composable arrows $f, g, h$, the equation $(h \circ g) \circ f = h \circ (g \circ f)$ holds.
    \item \textbf{Identity:} For each object $A$, there exists an identity morphism $\mathrm{id}_A: A \rightarrow A$ such that for any morphisms $f: A \rightarrow B$ and $g: B \rightarrow A$, we have $\mathrm{id}_B \circ f = f$ and $g \circ \mathrm{id}_A = g$.
\end{enumerate}

\section{Examples of Categories}
\begin{itemize}
    \item \textbf{Set:} Objects are sets, and morphisms are functions between sets.
    \item \textbf{Grp:} Objects are groups, and morphisms are group homomorphisms.
    \item \textbf{Top:} Objects are topological spaces, and morphisms are continuous maps.
    \item \textbf{Vect:} Objects are vector spaces, and morphisms are linear transformations.
\end{itemize}

\chapter{Slice Categories \cite{Jacobs}}
\section{Sets/I and Sets$(\rightarrow)$ taken from Jacobs}
Consider a family of sets as a function $\varphi : X \rightarrow I$. We often describe a family of sets as a function $\varphi : X \rightarrow I$ and say that $X$ is a family over $I$ and that $\varphi$ displays the family $(X_i)$. In order to emphasize that we think of such a map $\varphi$ as a family, we often write it vertically as $\left(\begin{smallmatrix} X \\ \downarrow \\ I \end{smallmatrix}\right)$. A constant family is one of the form $\left(\begin{smallmatrix} I \times X \\ \downarrow \\ I \end{smallmatrix}\right)$, where $\pi$ is the Cartesian product projection; often it is written simply as $I \times X$. Notice that all fibers of this constant family are (isomorphic to) $X$.

Such families $\left(\begin{smallmatrix} X \\ \downarrow \\ I \end{smallmatrix}\right)$ of sets give rise to two categories: the slice category $\mathbf{Sets}/I$ and the arrow category $\mathbf{Sets}(\rightarrow)$. The objects of $\mathbf{Sets}/I$ are the $I$-indexed families, for a fixed set $I$; the objects of $\mathbf{Sets}(\rightarrow)$ are all the $I$-indexed families, for all possible $I$.

\subsection{Sets/I}
Objects: families $\left(\begin{smallmatrix} X \\ \downarrow \\ I \end{smallmatrix}\right)$.

Morphisms: $\left(\begin{smallmatrix} X \\ \downarrow \\ I \end{smallmatrix}\right) \xrightarrow{f} \left(\begin{smallmatrix} Y \\ \downarrow \\ I \end{smallmatrix}\right)$ are functions $f : X \rightarrow Y$ making the following diagram commute.
\[
\begin{tikzcd}
X \arrow[rr, "f"] \arrow[dr, "\varphi"'] & & Y \arrow[dl, "\psi"] \\
 & I &
\end{tikzcd}
\]

\subsection{Sets$(\rightarrow)$}
Objects: families $\left(\begin{smallmatrix} X \\ \downarrow \\ I \end{smallmatrix}\right)$, for arbitrary sets $I$.

Morphisms: $\left(\begin{smallmatrix} X \\ \downarrow \\ I \end{smallmatrix}\right) \xrightarrow{(u, f)} \left(\begin{smallmatrix} Y \\ \downarrow \\ J \end{smallmatrix}\right)$ are pairs of functions $u : I \rightarrow J$ and $f : X \rightarrow Y$ for which the following diagram commutes.
\[
\begin{tikzcd}
 X \arrow[r, "f"] \arrow[d, "\varphi"] &
    Y \arrow[d, "\psi" ] \\
  I \arrow[r,  "u"] &
  J
\end{tikzcd}
\]

\noindent Hence, objects in the arrow category $\mathbf{Sets}(\rightarrow)$ involve an extra function $u$ between the index sets. Notice that one can now view $f$ as a collection of functions $f_i : X_i \rightarrow Y_{u(i)}$, since for $x \in \varphi^{-1}(i)$, $f(x)$ lands in $\psi^{-1}(u(i))$. Composition and identities in $\mathbf{Sets}(\rightarrow)$ are component-wise inherited from $\mathbf{Sets}$.

We further remark that there is a codomain functor $\text{cod} : \mathbf{Sets}(\rightarrow) \rightarrow \mathbf{Sets}$; it maps $\left(\begin{smallmatrix} X \\ \downarrow \\ I \end{smallmatrix}\right) \mapsto I$ and $\left(\begin{smallmatrix} (u, f) \end{smallmatrix}\right) \mapsto u$.

Also, for each $I$, there is a (non-full) inclusion functor $\mathbf{Sets}/I \rightarrow \mathbf{Sets}(\rightarrow)$.

\chapter{Monics, Epics, and Isomorphisms}
An arrow $f: A \rightarrow B$ in a category $\mathcal{C}$ is:
\begin{itemize}
    \item \emph{Monic} (monomorphism) if for all arrows $g, h: C \rightarrow A$, $f \circ g = f \circ h$ implies $g = h$.
    \item \emph{Epic} (epimorphism) if for all arrows $g, h: B \rightarrow C$, $g \circ f = h \circ f$ implies $g = h$.
    \item \emph{Isomorphic} (isomorphism) if there exists an arrow $g: B \rightarrow A$ such that $g \circ f = \mathrm{id}_A$ and $f \circ g = \mathrm{id}_B$.
\end{itemize}

\section{Examples in Set}
\begin{itemize}
    \item \textbf{Monic:} Injective functions.
    \item \textbf{Epic:} Surjective functions.
    \item \textbf{Isomorphisms:} Bijective functions.
\end{itemize}

\section{Examples in Grp}
\begin{itemize}
    \item \textbf{Monic:} Injective group homomorphisms.
    \item \textbf{Epic:} Surjective group homomorphisms.
    \item \textbf{Isomorphisms:} Bijective group homomorphisms (isomorphisms of groups).
\end{itemize}

\chapter{Diagrams, Cones, Cocones, Limits, Colimits}
In category theory, diagrams provide a structured way to visualize and understand the relationships between objects and morphisms. Cones and cocones are structures that relate diagrams to limits and colimits, respectively.

\section{Diagrams}
A \emph{diagram} in a category $\mathcal{C}$ is a functor $D: J \rightarrow \mathcal{C}$ where $J$ is an indexing category. Objects of $J$ are mapped to objects of $\mathcal{C}$, and morphisms of $J$ are mapped to morphisms in $\mathcal{C}$.

\subsection{Examples of Diagrams}
Diagrams can range from simple single-object diagrams to complex networks of objects and morphisms, such as sequences, commutative squares, or more intricate structures. They serve as a framework for discussing the relationships between objects and morphisms.
\[
\begin{tikzcd}
A \arrow[r, "f"] \arrow[d, "g"'] & B \arrow[d, "h"] \\
C \arrow[r, "k"] & D
\end{tikzcd}
\]

\section{Cones}
A \emph{cone} over a diagram $D: J \rightarrow \mathcal{C}$ consists of an object $L$ and a family of morphisms $\lambda_X: L \rightarrow D(X)$ for each object $X$ in $J$ such that for every morphism $f: X \rightarrow Y$ in $J$, the following diagram commutes:
\[
\begin{tikzcd}
& L \arrow[dl, "\lambda_X"'] \arrow[dr, "\lambda_Y"] \\
D(X) \arrow[rr, "D(f)"'] & & D(Y)
\end{tikzcd}
\]

\section{Limits}
``Let $ F : J \rightarrow \mathcal{C} $ be a diagram of shape $ J $ in a category $ \mathcal{C} $. A \textit{cone} to $ F $ is an object $ N $ of $ \mathcal{C} $ together with a family $ \psi_X : N \rightarrow F(X) $ of morphisms indexed by the objects $ X $ of $ J $, such that for every morphism $ f : X \rightarrow Y $ in $ J $, we have $ F(f) \circ \psi_X = \psi_Y $.

A \textit{limit} of the diagram $ F : J \rightarrow \mathcal{C} $ is a cone $ (L, \phi) $ to $ F $ such that for every cone $ (N, \psi) $ to $ F $, there exists a unique morphism $ u : N \rightarrow L $ such that $ \phi_X \circ u = \psi_X $ for all $ X $ in $ J $." \cite{Wikipedia}

\begin{center}
\begin{tikzcd}
 & N \arrow[ddl, "\psi_X"'] \arrow[ddr, "\psi_Y"] \arrow[dd, dashed, "u" description] \\
 & & \\
F(X) \arrow[rr, bend right, "F(f)"] & L \arrow[l, "\phi_X"'] \arrow[r, "\phi_Y"] & F(Y)
\end{tikzcd}
\end{center}

``One says that the cone $ (N, \psi) $ factors through the cone $ (L, \phi) $ with the unique factorization $ u $. The morphism $ u $ is sometimes called the \textit{mediating morphism}." \cite{Wikipedia}

\subsection{Examples of Limits}
\begin{itemize}
    \item \textbf{Products:} The product of two objects $A$ and $B$ in a category $\mathcal{C}$ is a limit of the diagram consisting of $A$ and $B$ with no morphisms between them.
    \item \textbf{Equalizers:} An equalizer of two parallel morphisms $f, g: A \rightarrow B$ is a limit of the diagram formed by $A$ and $B$ with two parallel arrows.
\end{itemize}

\section{Cocones}
A \emph{cocone} under a diagram $D: J \rightarrow \mathcal{C}$ consists of an object $L$ and a family of morphisms $\lambda_X: D(X) \rightarrow L$ for each object $X$ in $J$ such that for every morphism $f: X \rightarrow Y$ in $J$, the following diagram commutes:
\[
\begin{tikzcd}
D(X) \arrow[rr, "D(f)"] \arrow[dr, "\lambda_X"'] & & D(Y) \arrow[dl, "\lambda_Y"] \\
& L &
\end{tikzcd}
\]

\section{Colimits}
A \emph{colimit} of a diagram $D: J \rightarrow \mathcal{C}$ is a cocone $(L, \lambda)$ under $D$ that is universal among all such cocones. This means that for any other cocone $(L', \lambda')$ under $D$, there exists a unique morphism $u: L \rightarrow L'$ such that $\lambda'_X \circ u = \lambda_X$ for all $X$ in $J$.
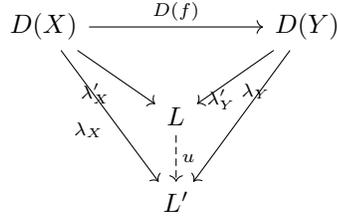
\[
\begin{tikzcd}
D(X) \arrow[rr, "D(f)"] \arrow[dr, "\lambda'_X"'] \arrow[ddr, "\lambda_X"'] & & D(Y) \arrow[dl, "\lambda_Y"] \arrow[ddl, "\lambda'_Y"'] \\
& L \arrow[d, "u", dashed] \\
& L' &
\end{tikzcd}
\]

\subsection{Examples of Colimits}
\begin{itemize}
    \item \textbf{Coproducts:} The coproduct of two objects $A$ and $B$ in a category $\mathcal{C}$ is a colimit of the diagram consisting of $A$ and $B$ with no morphisms between them.
    \item \textbf{Coequalizers:} A coequalizer of two parallel morphisms $f, g: A \rightarrow B$ is a colimit of the diagram formed by $A$ and $B$ with two parallel arrows.
\end{itemize}

\section{Equalizers and Coequalizers}
\subsection{Equalizers}
An \textit{equalizer} of two parallel arrows $f, g: A \rightarrow B$ is an object $E$ together with a morphism $e: E \rightarrow A$ such that $f \circ e = g \circ e$ and for any object $Z$ with a morphism $z: Z \rightarrow A$ such that $f \circ z = g \circ z$, there exists a unique morphism $u: Z \rightarrow E$ such that $e \circ u = z$. This is depicted as follows:
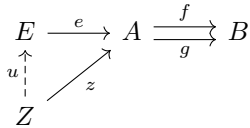
\[
\begin{tikzcd}
E \arrow[r, "e"] & A \arrow[r, shift left, "f"] \arrow[r, shift right, "g"'] & B \\
Z \arrow[u, dashed, "u"] \arrow[ur, "z"']
\end{tikzcd}
\]

\subsection{Coequalizers}
A \textit{coequalizer} of two parallel arrows $f, g: A \rightarrow B$ in a category $\mathcal{C}$ is an object $Q$ together with a morphism $q: B \rightarrow Q$, universal with the property $q \circ f = q \circ g$, as in the following diagram:
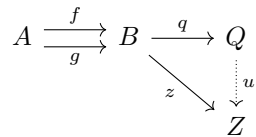
\cite{Awodey}
\[
\begin{tikzcd}
A \arrow[r,shift left=.75ex,"f"]
  \arrow[r,shift right=.75ex,swap,"g"]
&
B \arrow[r,"q"] \arrow[dr,swap,"z"]
&
Q \arrow[d,densely dotted,"u"]
\\
& & Z
\end{tikzcd}
\]
Given any $Z$ and $z: B \rightarrow Z$, if $z \circ f = z \circ g$, then there exists a unique $u: Q \rightarrow Z$ such that $u \circ q = z$ \cite{Awodey}.

\section{Products and coproducts}
\subsection{Products}
A \textit{product} of two objects $A$ and $B$ in a category $\mathcal{C}$ is an object $P$ together with two morphisms $\pi_1: P \rightarrow A$ and $\pi_2: P \rightarrow B$ such that for any object $X$ with morphisms $f: X \rightarrow A$ and $g: X \rightarrow B$, there exists a unique morphism $u: X \rightarrow P$ such that $\pi_1 \circ u = f$ and $\pi_2 \circ u = g$. This can be depicted as:
\[
\begin{tikzcd}
 & X \arrow[dl, "f"'] \arrow[dr, "g"] \arrow[d, dashed, "u" description] \\
A & P \arrow[l, "\pi_1"] \arrow[r, "\pi_2"'] & B
\end{tikzcd}
\]

\subsection{Coproducts}
A \textit{coproduct} of two objects $A$ and $B$ in a category $\mathcal{C}$ is an object $C$ together with two morphisms $\iota_1: A \rightarrow C$ and $\iota_2: B \rightarrow C$ such that for any object $X$ with morphisms $f: A \rightarrow X$ and $g: B \rightarrow X$, there exists a unique morphism $u: C \rightarrow X$ such that $u \circ \iota_1 = f$ and $u \circ \iota_2 = g$. This can be depicted as:
\[
\begin{tikzcd}
A \arrow[r, "\iota_1"] \arrow[dr, "f"'] & C \arrow[d, dashed, "u" description] & B \arrow[l, "\iota_2"'] \arrow[dl, "g"] \\
 & X &
\end{tikzcd}
\]

\section{Pushouts and Pullbacks}
\subsection{Pullbacks}
A \textit{pullback} (also known as a fiber product) of two morphisms $f: X \rightarrow Z$ and $g: Y \rightarrow Z$ in a category $\mathcal{C}$ is an object $P$ together with two morphisms $p_1: P \rightarrow X$ and $p_2: P \rightarrow Y$ such that the following diagram commutes:
\[
\begin{tikzcd}
P \arrow[r, "p_2"] \arrow[d, "p_1"'] & Y \arrow[d, "g"] \\
X \arrow[r, "f"'] & Z
\end{tikzcd}
\]
Moreover, $P$ must be universal with respect to this property, meaning that for any other object $Q$ with morphisms $q_1: Q \rightarrow X$ and $q_2: Q \rightarrow Y$ making the diagram commute, there exists a unique morphism $u: Q \rightarrow P$ such that $p_1 \circ u = q_1$ and $p_2 \circ u = q_2$. This situation is illustrated in the following commutative diagram:
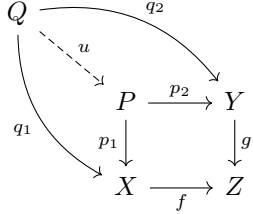
\[
\begin{tikzcd}
  Q \arrow[rrd,bend left,"q_2"]
    \arrow[ddr,bend right,swap,"q_1"]
    \arrow[dr,dashed,"u"] \\
  & P \arrow[d,"p_1"'] \arrow[r,"p_2"] & Y \arrow[d,"g"]  \\
  & X \arrow[r,swap,"f"]  & Z
\end{tikzcd}
\]
Explicitly, a pullback of the morphisms $f$ and $g$ consists of an object $P$ and two morphisms $p_1: P \rightarrow X$ and $p_2: P \rightarrow Y$ for which the diagram
\[
\begin{tikzcd}
P \arrow[r, "p_2"] \arrow[d, "p_1"'] & Y \arrow[d, "g"] \\
X \arrow[r, "f"'] & Z
\end{tikzcd}
\]
commutes. Moreover, the pullback $(P, p_1, p_2)$ must be universal with respect to this diagram. That is, for any other such triple $(Q, q_1, q_2)$ where $q_1: Q \rightarrow X$ and $q_2: Q \rightarrow Y$ are morphisms with $f \circ q_1 = g \circ q_2$, there must exist a unique $u: Q \rightarrow P$ such that
\[
p_1 \circ u = q_1, \quad p_2 \circ u = q_2.
\]
This situation is illustrated in the following commutative diagram:
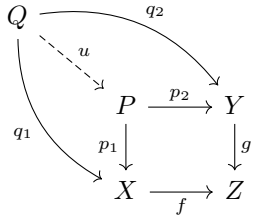
\[
\begin{tikzcd}
  Q \arrow[rrd,bend left,"q_2"]
    \arrow[ddr,bend right,swap,"q_1"]
    \arrow[dr,dashed,"u"] \\
  & P \arrow[d,"p_1"'] \arrow[r,"p_2"] & Y \arrow[d,"g"]  \\
  & X \arrow[r,swap,"f"]  & Z
\end{tikzcd}
\]

\subsection{Pushouts}
A \textit{pushout} (also known as a cofiber product) of two morphisms $f: Z \rightarrow X$ and $g: Z \rightarrow Y$ in a category $\mathcal{C}$ is an object $P$ together with two morphisms $\iota_1: X \rightarrow P$ and $\iota_2: Y \rightarrow P$ such that the following diagram commutes:
\[
\begin{tikzcd}
Z \arrow[r, "g"] \arrow[d, "f"'] & Y \arrow[d, "\iota_2"] \\
X \arrow[r, "\iota_1"'] & P
\end{tikzcd}
\]
Moreover, $P$ must be universal with respect to this property, meaning that for any other object $Q$ with morphisms $q_1: X \rightarrow Q$ and $q_2: Y \rightarrow Q$ making the diagram commute, there exists a unique morphism $u: P \rightarrow Q$ such that $u \circ \iota_1 = q_1$ and $u \circ \iota_2 = q_2$.
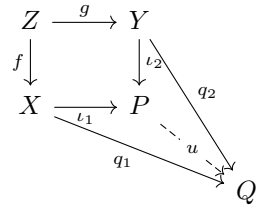
\[
\begin{tikzcd}
Z \arrow[r, "g"] \arrow[d, "f"'] & Y \arrow[d, "\iota_2"] \arrow[ddr, "q_2"] \\
X \arrow[r, "\iota_1"'] \arrow[drr, "q_1"'] & P \arrow[dr, dashed, "u" description] \\
& & Q
\end{tikzcd}
\]

\chapter{Initial and Terminal Objects \cite{Pettigrew}}
Consider the empty diagram in the category $\mathcal{C}$. The cones and cocones over this diagram are simply the objects of $\mathcal{C}$. Therefore:
\begin{itemize}
    \item If the empty diagram has a limit, it is an object $1$ such that, for every object $A$ in $\mathcal{C}$, there is a unique morphism $1_A$ : $A \rightarrow 1$.
    \item If the empty diagram has a colimit, it is an object $0$ such that, for every object $A$ in $\mathcal{C}$, there is a unique morphism $0_A$ : $0 \rightarrow A$.
\end{itemize}

\section{Definitions}
\begin{itemize}
    \item \textbf{Initial Object:} A limit of the empty diagram (if it exists) is called a terminal object of $\mathcal{C}$.
    \item \textbf{Terminal Object:} A colimit of the empty diagram (if it exists) is called an initial object of $\mathcal{C}$.
\end{itemize}

By Propositions 5.1.3 and 5.2.3, a terminal or initial object is unique up to isomorphism. By Axiom 2, every topos has initial and terminal objects.

\subsection{Examples}
\begin{itemize}
    \item In \textbf{Set}, the empty set $\emptyset$ is the only initial object, and any singleton set $\{a\}$ is a terminal object.
    \item In \textbf{Grp}, the trivial group $\{e_G\}$ is both initial and terminal. Such objects are called zero objects.
    \item In a category based on a poset, any minimum element is an initial object (if it exists), and any maximum element is a terminal object (if it exists).
    \item If $\mathcal{C}$ is a category, then $Id_X : X \rightarrow X$ is a terminal object of $\mathcal{C}/X$. If $\mathcal{C}$ has an initial object $0$, then $0_X : 0 \rightarrow X$ is an initial object of $\mathcal{C}/X$.
\end{itemize}

\section{Properties}
\begin{itemize}
    \item If $1$ is a terminal object and $f : 1 \rightarrow A$, then the arrow is monic.
    \item If $1$ is a terminal object, then $1 \times A \cong A \cong A \times 1$.
\end{itemize}

\chapter{Members of Objects \cite{Pettigrew}}
In traditional set theory, the fundamental notion is the membership relation. In category theory, the analogous concept involves arrows rather than elements. In a category $\mathcal{C}$ with a terminal object $1$, the members of an object $A$ are the morphisms from $1$ to $A$.

\section{Membership}
In the context of category theory, the action of picking out a member of a set $A$ can be understood as a function from a singleton set into $A$. This is because any singleton is a terminal object and any terminal object is a singleton in \textbf{Set}.

\textbf{Member of:} If $\mathcal{C}$ is a category with a terminal object $1$ and $A$ is an object of $\mathcal{C}$, then a member of $A$ is an arrow $x$ : $1 \rightarrow A$.

This definition implies that there cannot be 'membership chains' as members are arrows, not objects, and arrows cannot have members.

\begin{itemize}
    \item $1$ has exactly one element.
    \item If $0$ has an element, then $0 \cong 1$.
\end{itemize}

\section{Injective and Surjective}
We can also define injective and surjective arrows in this context.

\textbf{Injective and Surjective:} Given $f: A \rightarrow B$, we say that
\begin{itemize}
    \item The arrow is \textbf{injective} if, for all $x,y$ : $1 \rightarrow A$, if $fx = fy$ then $x = y$.
    \item The arrow is \textbf{surjective} if, for all $y : 1 \rightarrow B$, there is $x: 1 \rightarrow A$ such that $fx = y$.
\end{itemize}

In \textbf{Set}, the injective arrows are the monics, and the surjective arrows are the epics. However, this is not necessarily true in all categories with terminal objects.

\chapter{Exponential Objects \cite{Pettigrew}}
Axiom 2 ensures that all toposes have analogues of the addition and multiplication operations on sets—they are coproducts and products, respectively. However, it does not guarantee analogues of power sets $P(A)$ or function spaces $B^A = \{f: A \rightarrow B\}$. To address this, we introduce an axiom guaranteeing these objects.

Consider a (set-theoretical) function $f: A \times C \rightarrow B$. For every element $c \in C$, the function $f_c: a \mapsto f(a, c)$ is a function from $A$ to $B$. The function $\hat{f}: c \mapsto f_c$ maps $C$ into $B^A$. There is also a function $\mathrm{ev}: A \times B^A \rightarrow B$ that evaluates a function from $A$ to $B$ at a value in $A$.

Therefore:
\[
ev(a, f_c) = f_c(a) = f(a,c)
\]

\textbf{Exponential:} Suppose $\mathcal{C}$ is a category with products. For any objects $A$ and $B$, an exponential of $A$ and $B$ consists of
\begin{itemize}
    \item An object $B^A$ of category $\mathcal{C}$
    \item An arrow $ev : A \times B^A \rightarrow B$ of category $\mathcal{C}$
\end{itemize}
such that for any arrow $f: A \times C \rightarrow B$, there is an arrow $\hat{f} : C \rightarrow B^A$ making the following diagram commute:
\[
\begin{tikzcd}
A \times B^A \arrow[r, "\mathrm{ev}"] & B \\
A \times C \arrow[u, "\mathrm{Id}_A \times \hat{f}"] \arrow[ur, "f"']
\end{tikzcd}
\]

Exponentials of $A$ and $B$ are unique up to isomorphism, and anything isomorphic to an exponential of $A$ and $B$ is itself an exponential of $A$ and $B$.

We can check the validity of the exponential $B^A$ by considering whether it satisfies certain basic conditions stated in terms of the membership relation. In set theory, an exponential object should have the following property:
\[
f \text{ is a member of } B^A \text{ iff } f: A \rightarrow B
\]

Although this is not exactly the result we get in category theory, we do obtain a close analogue: there is a one-to-one correspondence between arrows $f : A \rightarrow B$ and arrows $g : 1 \rightarrow B^A$.

\textbf{Name of an arrow:} If $A \rightarrow B$ is an arrow, the name of $f$ (written $\ulcorner f \urcorner$) is the arrow $1 \rightarrow B^A$ such that the following diagram commutes:
\[
\begin{tikzcd}
A \times B^A \arrow[r, "\mathrm{ev}"] & B \\
A \times 1 \arrow[u, "\mathrm{Id}_A \times \ulcorner f \urcorner"] \arrow[ur, "f"']
\end{tikzcd}
\]

There is a bijection $f \mapsto \ulcorner f \urcorner$ between the set of arrows $f: A \rightarrow B$ and the set of arrows $\ulcorner f \urcorner : 1 \rightarrow B^A$.

With the definition of an exponential in hand, we can introduce the third axiom of a topos.

A topos has exponentials for every pair of objects.

\textbf{Cartesian Closed Category:} A category with limits for all finite diagrams and exponentials for all pairs of objects is called a Cartesian closed category.

\cite{Pettigrew}.

\chapter{Subobjects and Their Classifiers \cite{Pettigrew}}
\section{Subobjects}
In category theory, the analogue of a subset of a set $A$ is called a subobject of $A$. A subobject is not another object but an arrow, specifically a monic arrow from an object $S$ into $A$.

\textbf{Part or Subobject:} Suppose $A$ is an object. Then a subobject of $A$ is a monic arrow $S \rightarrow A$.

\textbf{Inclusion:} Suppose $S \rightarrow A$ and $T \rightarrow A$ are subobjects of $A$. We say that $S$ is included in $T$ (written $S \subseteq T$) if there is a morphism $S \rightarrow T$ such that the following diagram commutes:
\[
\begin{tikzcd}
S \arrow[dr] \arrow[r] & T \arrow[d] \\
& A
\end{tikzcd}
\]

\section{Subobject Classifiers}
Next, we introduce an object $\Omega$ in a topos whose members act as truth values and associate with each subobject $S \rightarrow A$ a characteristic function $A \rightarrow \Omega$.

\subsection{Motivation}
Consider sets where $\Omega = \{\text{true, false}\}$. The characteristic function of a subset $S \subseteq A$ can be seen as the function $\chi_S$ mapping $A$ to $\{\text{true, false}\}$ such that the inverse image of $\{\text{true}\}$ under $\chi_S$ is $S$.

\subsection{Inverse Images of Subobjects}
In \textbf{Set}, a pullback for the diagram
\[
\begin{tikzcd}[column sep=scriptsize]
B \arrow[dr , "f"]
& & S \arrow[dl, "i"] \\
& A
\end{tikzcd}
\]
is the set
\[
\{(b, s) \mid f(b) = i(s)\}
\]

If $S \rightarrow A$ and $i$ is the inclusion map $i: s \mapsto s$ for all $s \in S$, then a pullback is
\[
\{(b, s) \in B \times S \mid f(b) = i(s) = s\} \cong \{b \in B \mid f(b) \in S\} = f^{-1}(S)
\]

\textbf{Inverse Image:} If $S \rightarrow A$ and $B \rightarrow A$, then the pullback of the diagram
\[
\begin{tikzcd}
B \arrow[dr , "f"]
& & S \arrow[dl, "i"] \\
& A
\end{tikzcd}
\]
is called an inverse image of $S \rightarrow A$ under $f$ (written $f^{-1}(S)$).

\subsection{Definition}
We define the characteristic function $\chi_i$ of a subobject $S \rightarrow A$ to be the unique function such that $S$ is an inverse image of $1 \rightarrow \Omega$ under $\chi_i$.

\textbf{Subobject Classifier:} Suppose $\mathcal{C}$ is a category with a terminal object $1$. A subobject classifier in $\mathcal{C}$ consists of:
\begin{itemize}
    \item An object $\Omega$ of category $\mathcal{C}$.
    \item An arrow $1 \rightarrow \Omega$
\end{itemize}
such that for any object $A$ and subobject $S \rightarrow A$, there is a unique arrow $\chi_S: A \rightarrow \Omega$ such that:
\[
\begin{tikzcd}
S \arrow[r] \arrow[d] & 1 \arrow[d, "1"] \\
A \arrow[r, "\chi_i"] & \Omega
\end{tikzcd}
\]
is a pullback square.

Subobject classifiers are unique up to isomorphism, and anything isomorphic to a subobject classifier is itself a subobject classifier.

The final axiom of toposes can now be stated:

A topos has a subobject classifier.

\textbf{False:} The arrow $1 \rightarrow \Omega$ is the characteristic function of the subobject $0 \rightarrow 1$.

\cite{Pettigrew}.

\chapter{The Definition of a Topos}
\section{The Definition}
Having outlined the necessary axioms, we define a topos:

\textbf{Topos:} A topos is a category with:
\begin{itemize}
    \item Limits and colimits for all finite diagrams.
    \item Exponentials for every pair of objects.
    \item A subobject classifier.
\end{itemize}

\section{Examples}
\begin{itemize}
    \item \textbf{Set} is a topos.
    \item \textbf{Grp} is not a topos.
    \item \textbf{FinSet} is a topos.
\end{itemize}

\section{Fundamental Theorem of Toposes}
\textit{``If $\mathcal{E}$ is a topos and X is an object in $\mathcal{E}$, then $\mathcal{E}$/X is a topos."}
\cite{Pettigrew}.

\chapter{Algebra of Subobjects}

In category theory, understanding subobjects and their algebra is crucial for grasping the structure and properties of categories. Here, we explain subobject lattices and their relationship to Boolean and Heyting algebras.

\subsection{Subobject Lattices}

A \textbf{subobject} of an object \( A \) in a category \( \mathcal{C} \) is defined as an equivalence class of monomorphisms \( m: S \hookrightarrow A \). Two monomorphisms \( m: S \hookrightarrow A \) and \( m': S' \hookrightarrow A \) are considered equivalent if there exists an isomorphism \( f: S \rightarrow S' \) such that \( m = m' \circ f \).

The collection of all subobjects of \( A \), denoted as \(\text{Sub}(A)\), forms a \textbf{partially ordered set} (poset) under inclusion. Specifically, for subobjects \( [m: S \hookrightarrow A] \) and \( [m': S' \hookrightarrow A \), we have \( [m] \leq [m'] \) if there exists a morphism \( f: S \rightarrow S' \) making the following diagram commute:
\[
\begin{tikzcd}
S \arrow[dr, "m"'] \arrow[rr, "f"] & & S' \arrow[dl, "m'"] \\
& A &
\end{tikzcd}
\]

This poset of subobjects, \(\text{Sub}(A)\), possesses additional structure, making it a \textbf{lattice}.

\subsection{Lattice Structure of Subobjects}

A \textbf{lattice} is a poset in which any two elements have a greatest lower bound (glb) or meet, and a least upper bound (lub) or join. For subobjects \( S \) and \( T \) of \( A \):

\begin{itemize}
    \item The \textbf{meet} \( S \wedge T \) is given by the intersection of \( S \) and \( T \) in \( A \).
    \item The \textbf{join} \( S \vee T \) is represented by the subobject generated by the union of \( S \) and \( T \).
\end{itemize}

Formally, these operations can be defined through pullbacks and pushouts in the category \( \mathcal{C} \).

\subsection{Heyting Algebra of Subobjects}

In a \textbf{Heyting category}, the poset \(\text{Sub}(A)\) is not just a lattice but a \textbf{Heyting algebra}. This means it supports an additional operation called \textbf{implication} \( \rightarrow \), which is characterized by the property that for subobjects \( S \) and \( T \) of \( A \), there is a largest subobject \( U \) such that \( S \wedge U \leq T \).

Heyting algebras generalize Boolean algebras and are essential in the internal logic of a topos, particularly in intuitionistic logic.

\subsection{Power Objects in Topos Theory}

A \textbf{topos} can be seen as a categorical generalization of set theory, incorporating both logical and geometrical aspects. In a topos, every subobject poset \(\text{Sub}(A)\) is a Heyting algebra, and there exists a \textbf{power object} \( \mathcal{P}(A) \), which internalizes the notion of the power set in set theory.

The power object \( \mathcal{P}(A) \) of an object \( A \) is equipped with a monomorphism \( \in_A : A \times \mathcal{P}(A) \rightarrow \Omega \), where \( \Omega \) is the subobject classifier of the topos. This morphism satisfies a universal property analogous to the characteristic function of a subset in set theory.

\subsection{Boolean Topos}

In a \textbf{Boolean topos}, the internal Heyting algebra \(\text{Sub}(A)\) for any object \( A \) is a Boolean algebra. This implies that every subobject has a complement, and the internal logic of the topos aligns with classical Boolean logic.

\subsection{Internal Logic and Lawvere-Tierney Topology}

The internal logic of a topos can be intuitionistic or classical, depending on whether it forms a Heyting algebra or a Boolean algebra. The \textbf{Lawvere-Tierney topology} provides a framework to study these logical structures within a topos, defining modalities that extend the logical operations within the category.

\section*{References}

\begin{itemize}
    \item Pettigrew, R. \emph{An Introduction to Toposes}. [Unpublished Notes].
    \item Awodey, S. (2010). \emph{Category Theory}. Oxford Logic Guides 52, Oxford University Press.
    \item \href{https://ncatlab.org/nlab/show/Heyting+algebra#to_toposes}{nLab: Heyting Algebra}
\end{itemize}

\chapter{Kinds of Topos}
\section{Non-degenerate Toposes}
A topos is non-degenerate if it has more than one object and more than one morphism.

\section{Well-pointed Toposes}
A topos is well-pointed if the only arrow $1 \rightarrow 1$ is the identity arrow.

\section{Bivalent Toposes}
A topos is bivalent if it has exactly two objects and two morphisms.

\section{Boolean Toposes}
A topos is Boolean if its subobject classifier is a Boolean algebra.

\section{Natural Number Objects}
A topos has a natural number object if there is an object $N$ and arrows $0: 1 \rightarrow N$ and $s: N \rightarrow N$ such that for any object $A$ and arrows $f: 1 \rightarrow A$ and $g: A \rightarrow A$, there is a unique arrow $h: N \rightarrow A$ such that $h \circ 0 = f$ and $h \circ s = g \circ h$.

\cite{Pettigrew}.

\chapter{Functors}
\section{The Definition of a Functor}
A functor $F: \mathcal{C} \rightarrow \mathcal{D}$ consists of:
\begin{itemize}
    \item A function $F: \mathrm{Ob}(\mathcal{C}) \rightarrow \mathrm{Ob}(\mathcal{D})$
    \item A function $F: \mathrm{Ar}(\mathcal{C}) \rightarrow \mathrm{Ar}(\mathcal{D})$
\end{itemize}
such that the following properties hold:
\begin{itemize}
    \item $F(\mathrm{id}_A) = \mathrm{id}_{F(A)}$ for all objects $A$ in $\mathcal{C}$.
    \item $F(g \circ f) = F(g) \circ F(f)$ for all arrows $f, g$ in $\mathcal{C}$.
\end{itemize}

\section{Examples of Functors}
\begin{itemize}
    \item The identity functor $\mathrm{Id}: \mathcal{C} \rightarrow \mathcal{C}$.
    \item The constant functor $C: \mathcal{C} \rightarrow \mathcal{D}$ that sends every object to a fixed object $D$ and every arrow to $\mathrm{id}_D$.
    \item The hom-functor $\mathrm{Hom}(A, -): \mathcal{C} \rightarrow \mathbf{Set}$ that sends every object $B$ to the set of arrows $A \rightarrow B$ and every arrow $f: B \rightarrow C$ to the function $\mathrm{Hom}(A, f): \mathrm{Hom}(A, B) \rightarrow \mathrm{Hom}(A, C)$ given by composition with $f$.
\end{itemize}

\section{Properties of Functors}
Functors preserve the structure of categories in various ways:
\begin{itemize}
    \item They preserve isomorphisms: If $f: A \rightarrow B$ is an isomorphism in $\mathcal{C}$, then $F(f)$ is an isomorphism in $\mathcal{D}$.
    \item They preserve commutative diagrams: If a diagram commutes in $\mathcal{C}$, its image under $F$ commutes in $\mathcal{D}$.
    \item They preserve limits and colimits: If $\mathcal{C}$ has a limit (or colimit) of a diagram $D$, then $\mathcal{D}$ has a limit (or colimit) of $F(D)$.
\end{itemize}

\chapter{Natural Transformations}
\section{The Definition of a Natural Transformation}
A natural transformation $\eta: F \Rightarrow G$ between two functors $F, G: \mathcal{C} \rightarrow \mathcal{D}$ consists of:
\begin{itemize}
    \item For each object $A$ in $\mathcal{C}$, an arrow $\eta_A: F(A) \rightarrow G(A)$ in $\mathcal{D}$
\end{itemize}
such that for every arrow $f: A \rightarrow B$ in $\mathcal{C}$, the following diagram commutes:
\[
\begin{tikzcd}
F(A) \arrow[r, "F(f)"] \arrow[d, "\eta_A"'] & F(B) \arrow[d, "\eta_B"] \\
G(A) \arrow[r, "G(f)"] & G(B)
\end{tikzcd}
\]

This commutative diagram ensures that natural transformations respect the structure of the categories involved, providing a way to compare functors in a coherent manner. 

\section{Examples of Natural Transformations}
\begin{itemize}
    \item The identity natural transformation $\mathrm{Id}: F \Rightarrow F$ assigns the identity morphism to each object, ensuring that $\eta_A = \mathrm{id}_{F(A)}$ for all $A$ in $\mathcal{C}$.
    \item The composition of two natural transformations $\eta: F \Rightarrow G$ and $\mu: G \Rightarrow H$ is a natural transformation $\mu \circ \eta: F \Rightarrow H$ defined by $(\mu \circ \eta)_A = \mu_A \circ \eta_A$ for each object $A$ in $\mathcal{C}$.
    \item For the hom-functor $\mathrm{Hom}(A, -): \mathcal{C} \rightarrow \mathbf{Set}$, a natural transformation between $\mathrm{Hom}(A, -)$ and another functor $\mathrm{Hom}(B, -)$ would be given by a function $\eta: \mathrm{Hom}(A, X) \rightarrow \mathrm{Hom}(B, X)$ for each $X$ in $\mathcal{C}$ that respects composition and identity.
\end{itemize}

\section{Properties of Natural Transformations}
Natural transformations have several key properties:
\begin{itemize}
    \item \textbf{Naturality:} The naturality condition, expressed by the commutative diagram above, ensures that the transformation is consistent with the action of the functors on morphisms.
    \item \textbf{Vertical and Horizontal Composition:} Natural transformations can be composed both vertically (as in the composition of $\eta$ and $\mu$) and horizontally, which involves transforming functors applied in sequence.
\end{itemize}

\cite{Pettigrew}.

\chapter{Adjoint Functors}
\section{Definition}
Two functors \(F: \mathcal{C} \rightarrow \mathcal{D}\) and \(G: \mathcal{D} \rightarrow \mathcal{C}\) are adjoint if there is a natural isomorphism
\[
\mathrm{Hom}_{\mathcal{D}}(F(A), B) \cong \mathrm{Hom}_{\mathcal{C}}(A, G(B))
\]
for all objects \(A\) in \(\mathcal{C}\) and \(B\) in \(\mathcal{D}\).

\subsection{Definition via Counit-Unit Adjunction \cite{Wikipedia}}
``A counit-unit adjunction between two categories \( \mathcal{C} \) and \( \mathcal{D} \) consists of two functors \( F : \mathcal{D} \rightarrow \mathcal{C} \) and \( G : \mathcal{C} \rightarrow \mathcal{D} \) and two natural transformations \( \varepsilon : FG \rightarrow 1_{\mathcal{C}} \) and \( \eta : 1_{\mathcal{D}} \rightarrow GF \), respectively called the counit and the unit of the adjunction, such that the compositions
\[
F \xrightarrow{F \eta} FGF \xrightarrow{\varepsilon F} F \quad \text{and} \quad G \xrightarrow{\eta G} GFG \xrightarrow{G \varepsilon} G
\]
are the identity transformations \( 1_F \) and \( 1_G \) on \( F \) and \( G \) respectively. In this situation, we say that \( F \) is left adjoint to \( G \) and \( G \) is right adjoint to \( F \), and may indicate this relationship by writing \( (\varepsilon, \eta) : F \dashv G \), or simply \( F \dashv G \).

In equation form, the above conditions on \( (\varepsilon, \eta) \) are the counit-unit equations:
\[
1_F = \varepsilon F \circ F \eta \quad \text{and} \quad 1_G = G \varepsilon \circ \eta G
\]
which mean that for each \( X \) in \( \mathcal{C} \) and each \( Y \) in \( \mathcal{D} \),
\[
1_{FY} = \varepsilon_{FY} \circ F(\eta_Y) \quad \text{and} \quad 1_{GX} = G(\varepsilon_X) \circ \eta_{GX}
\]

Note that \( 1_{\mathcal{C}} \) denotes the identity functor on the category \( \mathcal{C} \), \( 1_F \) denotes the identity natural transformation from the functor \( F \) to itself, and \( 1_{FY} \) denotes the identity morphism of the object \( FY \)."

\subsection{Definition via Universal Morphisms \cite{Wikipedia}}
``By definition, a functor \( F : \mathcal{D} \rightarrow \mathcal{C} \) is a left adjoint functor if for each object \( X \) in \( \mathcal{C} \) there exists a universal morphism from \( F \) to \( X \). Spelled out, this means that for each object \( X \) in \( \mathcal{C} \) there exists an object \( G(X) \) in \( \mathcal{D} \) and a morphism \( \varepsilon_X : F(G(X)) \rightarrow X \) such that for every object \( Y \) in \( \mathcal{D} \) and every morphism \( f : F(Y) \rightarrow X \) there exists a unique morphism \( g : Y \rightarrow G(X) \) with \( \varepsilon_X \circ F(g) = f \).

The latter equation is expressed by the following commutative diagram:
\begin{center}
\begin{tikzcd}
F(Y) \arrow[rr, dashed, "F(g)"] \arrow[dr, "f"'] & & F(G(X)) \arrow[dl, "\varepsilon_X"] \\
 & X &
\end{tikzcd}
\end{center}

In this situation, one can show that \( G \) can be turned into a functor \( G : \mathcal{C} \rightarrow \mathcal{D} \) in a unique way such that \( \varepsilon_X \circ F(G(f)) = f \circ \varepsilon_{X'} \) for all morphisms \( f : X' \rightarrow X \) in \( \mathcal{C} \); \( F \) is then called a left adjoint to \( G \).

Similarly, we may define right-adjoint functors. A functor \( G : \mathcal{C} \rightarrow \mathcal{D} \) is a right adjoint functor if for each object \( Y \) in \( \mathcal{D} \) there exists a universal morphism from \( Y \) to \( G \). Spelled out, this means that for each object \( Y \) in \( \mathcal{D} \) there exists an object \( F(Y) \) in \( \mathcal{C} \) and a morphism \( \eta_Y : Y \rightarrow G(F(Y)) \) such that for every object \( X \) in \( \mathcal{C} \) and every morphism \( g : Y \rightarrow G(X) \) there exists a unique morphism \( f : F(Y) \rightarrow X \) with \( G(f) \circ \eta_Y = g \)

The latter equation is expressed by the following commutative diagram:"
\begin{center}
\begin{tikzcd}
Y \arrow[rr, "\eta_Y"] \arrow[dr, "g"'] & & G(F(Y)) \arrow[dl, "G(f)", dashed] \\
 & G(X) &
\end{tikzcd}
\end{center}

\section{Examples Concerning Set}
\begin{itemize}
    \item The forgetful functor $U: \mathbf{Grp} \rightarrow \mathbf{Set}$ has a left adjoint $F: \mathbf{Set} \rightarrow \mathbf{Grp}$ that sends each set $S$ to the free group on $S$. This means there is a natural isomorphism $\mathrm{Hom}_{\mathbf{Grp}}(F(S), G) \cong \mathrm{Hom}_{\mathbf{Set}}(S, U(G))$ for any set $S$ and group $G$.
    \item The forgetful functor $U: \mathbf{Top} \rightarrow \mathbf{Set}$ has a left adjoint $F: \mathbf{Set} \rightarrow \mathbf{Top}$ that sends each set $S$ to the discrete topology on $S$. Here, $\mathrm{Hom}_{\mathbf{Top}}(F(S), T) \cong \mathrm{Hom}_{\mathbf{Set}}(S, U(T))$ for any set $S$ and topological space $T$.
\end{itemize}

\section{Examples Concerning Forgetful Functors}
Forgetful functors between various categories often have adjoints. For example:
\begin{itemize}
    \item The forgetful functor $U: \mathbf{Mon} \rightarrow \mathbf{Set}$ has a left adjoint $F: \mathbf{Set} \rightarrow \mathbf{Mon}$ that sends each set $S$ to the free monoid on $S$.
    \item The forgetful functor $U: \mathbf{Vect} \rightarrow \mathbf{Set}$ has a left adjoint $F: \mathbf{Set} \rightarrow \mathbf{Vect}$ that sends each set $S$ to the free vector space on $S$.
\end{itemize}
\cite{Pettigrew}.

\end{document}